\documentclass[12pt, draft]{article}
\usepackage{amsmath}
\usepackage{amssymb}
\usepackage{latexsym}
\usepackage{color}
\allowdisplaybreaks[4]

\newtheorem{thm}{\quad Theorem}[section]
\newtheorem{prp}[thm]{\quad Proposition}
\newtheorem{lmm}[thm]{\quad Lemma}

\newtheorem{rem}[thm]{\quad Remark}
\renewcommand{\v}[1]{\vert #1 \vert}
\newcommand{\V}[1]{\Vert #1 \Vert}

\newcommand{\N}{{\mathbb N}}

\newcommand{\R}{{\mathbb R}}

\title{The shortest lengths and the enclosure method for time dependent problems}
\author{
Mishio Kawashita\thanks{Partly supported by JSPS KAKENHI Grant Number JP19K03565.}
\thanks{kawasita@hiroshima-u.ac.jp}
\\
\and
Wakako Kawashita\thanks{Partly supported by JSPS KAKENHI Grant Number JP20K03684}
\thanks{wakawa@hiroshima-u.ac.jp}
}
\begin{document}
\maketitle
\begin{abstract}
In this article, the role of \lq\lq shortest length''  is discussed concerning the inverse problems with the time dependent enclosure method.  In the beginning, the inclusions embedded in a non-layered medium are handled. The argument for the above case reveals why the \lq\lq shortest length'' is important.  
Next, the case where the inclusions are embedded in a two-layered medium is explained in view of the previous works by Ikehata and the authors.   In addition, the regularity assumptions for the boundaries of the inclusions are relaxed.  
\end{abstract}
%

\vskip 1truepc
%
%
%
%
%
%
%
%
\setcounter{equation}{0}
\section{Introduction}

For mathematical models of nondestructive testing, 
the enclosure method is introduced by M. Ikehata \cite{Ikehata M 2000, Ikehata M 1999}  by using measurements
consisting of the pairing of both the voltage potential and electric current on the boundary. 
This problem is formulated by elliptic boundary value problems. 
In \cite{Ikehata M 2000, Ikehata M 1999}, Ikehata gives some procedure to find the convex hull $D'$ of cavities $D$, which are parts with no conductivity and embedded inside of media, by obtaining the support function $h(\omega) = \sup_{x \in D}x\cdot\omega$ ($\v{\omega} = 1$) 
of the cavities. Once the support function is known, the convex hull $D'$ is recovered as $D' = \bigcap_{\v{\omega} = 1}\{x  \vert x\cdot\omega < h(\omega) \}$. 
This is a new insight in the fields of the inverse problems governed by differential equations. 
\par
Now the enclosure method is investigated by many researchers, and one of the main streams of inverse problems is to handle the time dependent problems, which are called \lq\lq the time dependent enclosure method''. 
A survey for recent development of the time dependent enclosure method is given by Ikehata \cite{IR}.  
For inverse problems of differential equations including the 
enclosure method, 
see for examples \cite{Iksugaku, ISA3, Isozaki2, Nakamura-Potthast, Uhlmann} and 
earlier references given there.

In this article, we consider some prototype of inverse problems related with subsurface radar \cite{DGS} or ground probing radar \cite{BA}.  
Obstacles buried in the ground may consist of 
different media from background surrounding the obstacles. These are represented by 
inclusions which have different propagation speeds from the one of the background matrix media.
Mathematically, the simplest formulation is the case that the background consists of only one material and thus has no layers. 
A mathematical formulation in the context of the time dependent enclosure method is given by Ikehata \cite{IE, IE02}.

\par
For a formulation corresponding to buried obstacles, we take a two-layered medium with the flat boundary as a background. 
The obstacles are buried in the below of this transmission boundary. This is
the simplest case for the buried obstacles. This case is developed by \cite{transmission No1, transmission No2}.

\par

Before explaining the two-layer case, 
in section \ref{no layered case}, we consider the case of no layers.
Assuming that the boundary of the inclusions is smooth, 
Ikehata \cite{IE02} obtains the result itself essentially.
The formulations and the indicator functions given in \cite{IE02} are different from ours.  
In this article, we handle this problem by using the argument developed in \cite{transmission No1, transmission No2}. 
In the non-layer case, it is relatively simple to formulate the problems.
To obtain the main result (cf. Theorem \ref{CP case}),  
the most important purpose is to explain the shortest length or time for waves 
traveling from the obstacles $D$ to the observation set $B$ where we take measurement. 
The \lq\lq shortest length'' is appeared in showing 
asymptotic behavior of the resolvent deduced from the original time dependent problems. 
Note that in this article, solutions of some reduced problems corresponding to
the time dependent problems are also called \lq\lq resolvents'' though it seems to make an improper use of terminology. 

\par
We can give the integral kernel of resolvent for the reduced problems
of the usual wave equations
explicitly (cf. (\ref{Equation of v for CP}) and 
(\ref{integral kernel representation for CP case})). 
By this integral representation, we can easily know why the shortest length, which is the same as the shortest time in this case, can be found (cf. Proposition \ref{estimate of nabla_xv for CP}).

\par

Perhaps, a part of the main result (cf. Theorem \ref{CP case}) in section \ref{no layered case} is new. We know that the regularity assumption for the boundary of the inclusions can be reduced to a cone condition. Instead of that, we have to assume that the observation set $B$ is convex set. 
In \cite{transmission No1, transmission No2}, by using the regularity of the boundary, 
we show that the waves giving the shortest length return to the starting points. 
This property is also ensured under these relaxed assumptions.

\par

In section \ref{Inverse problems in two-layered medium}, 
we introduce our previous works \cite{transmission No1, transmission No2} on the inverse 
problems for obstacles buried in the below of the flat layer (i.e. 
the flat transmission boundary). 
For other approaches, see e.g. \cite{BKL, DEKPS, LLLL, LZ}.
In this case, waves are reflected and refracted by the layer. 
Snell's law gives the relation between incident waves and refracted waves, and determines the time it takes waves to go to a point $x$ in the lower side of the layer from $y$ in the upper side. 
This shortest time is called \lq\lq the optical distance'' between $x$ and $y$. 
From the point of view of section \ref{no layered case}, the optical distance can be found by asymptotic behavior of the resolvent (cf. Proposition \ref{estimate of nabla_xv for transmission case}).

\par

The main difference between the non-layer case and the two-layer case is the asymptotic behavior of the integral kernel of the resolvent (See Proposition \ref{Asymptotic behaviour of Phi_tau(x, y)}, which is shown in \cite{transmission No1, transmission No2}).  When we consider the layered medium, we have to take account of the total reflection phenomena.  If the speed of waves in the upper side is greater than that in the lower side, the evanescent waves may appear for the waves coming from the lower side. These waves propagate along the transmission boundary. 
The question is whether the evanescent waves may affect asymptotic behavior of the resolvent or not.

\par

In section \ref{how to show the case of less regularities}, 
we show that assumption for the regularities of the inclusion can be reduced as in section \ref{no layered case}. A similar argument works for the two-layer case. 
This is an extension of \cite{transmission No1, transmission No2}.

\par

In section \ref{Another choice of the indicator functions}, we consider another choice 
of the indicator function (cf. (\ref{another indicator function})). 
Eventually, this indicator function has the same work as that originally
given by (\ref{Indicator function for CP}).  
The original indicator (\ref{Indicator function for CP}) 
is easier to handle, but on the other hand another indicator function 
(\ref{another indicator function}) makes the problems clear.

\vskip1pc
\setcounter{equation}{0}
\section{Inverse problems for monotone inclusions in a non-layered medium}
\label{no layered case}

Let $\gamma(x) = (\gamma_{ij}(x))$ be a symmetric matrix whose components 
$\gamma_{ij}$ satisfy $\gamma_{ij} \in L^{\infty}(\R^3)$ $(i, j = 1, 2, 3)$.
Throughout this article, we assume that $\gamma$ is positive definite a.e. $x \in \R^3$, i.e. there exists a constant $c_0 > 0$ such that $\gamma(x) \geq c_0$ (which means that $\gamma(x)\xi\cdot\xi\geq c_0\vert\xi\vert^2$ ($\xi \in \R^3$)) a.e. $x \in \R^3$.
For $\gamma$, we put $L_{\gamma}u = {\rm div}(\gamma(x)\nabla_xu(x))$. 
\par

We fix $T > 0$ and consider the following wave equation: 
\begin{equation}
\left\{
\begin{array}{ll}
\displaystyle
(\partial_t^2-L_{\gamma}) u = 0 & \qquad\text{in}\, (0, T)\times\R^3,
\\
\displaystyle
u(0, x) = 0, \quad
\partial_tu(0, x) = f(x) &  \qquad\text{on}\,\R^3.
\end{array}
\right.
\label{equation with inclusion}
\end{equation}
In the case $\gamma = I_3$, where $I_3$ is the $3{\times}3$ unit matrix, (\ref{equation with inclusion}) is the Cauchy problem for the usual wave equation. 

\par

Note that for any $f \in L^2(\R^3)$, there exists a unique weak solution $u$
of (\ref{equation with inclusion}) satisfying $u \in L^2(0, T ; H^1(\R^3))$ with
$\partial_tu \in L^2(0, T ; L^2(\R^3))$ and $\partial_t^2u \in L^2(0, T ; (H^1(\R^3))')$
satisfying 
$$
\langle \partial_t^2u(t, \cdot), \phi \rangle 
+ \int_{\R^3}\gamma(x)\nabla_xu(t, x)\cdot\nabla_x\phi(x)dx = 0 \qquad \text{a.e. $t \in (0, T)$}$$
for all $\phi \in H^1(\R^3)$, and the initial condition in usual sense (cf. for example \cite{DL}). 
In what follows, we consider these weak solutions as the class of solutions.

\par
In (\ref{equation with inclusion}), taking $\gamma$ appropriately, we can formulate media 
containing inclusions 
$D$ standing for different media in the matrix material. Assume that $D$ is bounded, 
$\R^3\setminus{D}$ is connected, $D$ and $\gamma$ satisfy 
$\gamma(x) = I_3$ in $\R^3\setminus{D}$, and either of the following conditions: 
\par\noindent\hskip24pt
$(\text{M})_+$ : $\gamma$ satisfies $\inf_{x \in D}\gamma(x) > I_3$ 
\hskip12pt
or 
\hskip12pt
 $(\text{M})_-$ : $\gamma$ satisfies $\sup_{x \in D}\gamma(x) < I_3$.
\par\noindent
These a priori conditions are called \lq\lq monotonicity condition'' for $\gamma$.

In this setting, Ikehata \cite{IE, IE02} gives a formulation of one of prototype inverse problems for signal emission. For various cases with cavities or for taking measurements in different places, see also \cite{IE, IE02}. In this article, we only consider some typical case by using the approach handled in \cite{transmission No1, transmission No2} for the two-layer case.

\par

Take an open set $B$ with $\overline{B}\cap\overline{D} = \emptyset$, and put 
$f \in L^2(\R^3)$ satisfying the emission condition on $B$, which is given by
\begin{equation}
\left\{
\begin{array}{ll}
\text{$f \in L^2(\R^3)$ with ${\rm supp } f \subset \overline{B} $ and there exists a constant $c_1 > 0$ such that } 
 \\
\text{$f(x) \geq c_1$ ($x \in B$) or $-f(x) \geq c_1$ ($x \in B$).}
\end{array}
\right.
\label{emission condition}
\end{equation}
This condition (\ref{emission condition}) ensures that waves are emanated from $B$ exactly. 
Throughout this section, we always assume that (\ref{emission condition}) holds.
For this initial data $f$, we measure waves on $B$ from the initial time $0$ to $T$. Thus, the measurement is given by $u(t, x)$ for $0 \leq t \leq T$ and $x \in B$. Hence, in this setting,
the inverse problem to be considered is to find information of the inclusion $D$ from this measurement.

For this purpose, we introduce the indicator function $I_\tau$ defined by
\begin{equation}
I_\tau = \int_{\R^3}f(x)(w(x; \tau) - v(x; \tau))dx
\qquad(\tau \geq 1),
\label{Indicator function for CP}
\end{equation}
where $w(x; \tau)$ is defined by
\begin{equation}
w(x; \tau) = \int_{0}^{T}e^{-{\tau}t}u(t, x)dt
\label{def of w}
\end{equation}
and $v(x; \tau)$ is the solution of 
\begin{equation}
(L_{I_3}-\tau^2)v(x; \tau)+f(x) = 0 \qquad\text{in } \R^3.
\label{Equation of v for CP}
\end{equation}
Note that $I_\tau$ is obtained from the measurement $u(t, x)$ 
for $0 \leq t \leq T$ and $x \in B$.
\par

We put ${\rm dist}(D, B) = \inf_{x \in D, y \in B}\v{x - y}$. 
Since we only consider the case of inclusions, conditions for the boundary $\partial{D}$ of $D$ can
be relaxed. For $x \in \R^3$, $h > 0$, $\theta > 0$ 
and $n \in \R^3$ with $\v{n} = 1$, put
$C(x, n, h, \theta)  = \{y \in \R^3 \vert \v{y - x} \leq h, 
(y-x){\cdot}n \geq \v{y - x}\cos\theta \}$. Note that $C(x, n, h, \theta)$ is a cone with 
vertex $x$, direction $n$, radius $h$ and opening angle $\theta$.
We introduce the following condition (C):
\par\noindent
$$
{\rm (C)} \quad
\left\{
\begin{array}{ll}
\text{for any $x \in \partial{D}$, there exists a cone $C(x, n, h, \theta)$ with vertex $x$} 
\\
\text{satisfying $C(x, n, h, \theta)\setminus\{x\} \subset D$.}
\end{array}
\right.
$$
Condition (C) is a kind of cone conditions. 
Note that domains with $C^1$ boundary satisfy
condition (C). 

\par

From asymptotic behavior of $I_\tau$ as $\tau \to \infty$, 
we can obtain the following result which is essentially shown by 
Ikehata \cite{IE02} for the case that $\partial{D}$ is smooth:
\begin{thm}\label{CP case}
Suppose that $D$ and $B$ satisfy either $(a)$ or $(b)$: \\
$(a)$ $D$ has $C^1$ boundary $\partial{D}$ and $B$ satisfies condition $(C)$, \\
$(b)$ $D$ and $B$ satisfy condition $(C)$ and $B$ is convex. \\
Then, the indicator function defined by (\ref{Indicator function for CP}) satisfies 
\par\noindent
$(i)$ for $T < 2{\rm dist}(D, B)$, $\lim_{\tau \to \infty}e^{{\tau}T}I_\tau = 0$, 
\par\noindent
$(ii)$ for $T > 2{\rm dist}(D, B)$, $\lim_{\tau \to \infty}e^{{\tau}T}I_\tau = \mp\infty$ 
if $(\text{M})_\pm$ is assumed, respectively.
\par
Further, suppose $T > 2{\rm dist}(D, B)$ and either $(\text{M})_+$ or $(\text{M})_-$. Then
\begin{equation}
\lim_{\tau \to \infty}\frac{1}{\tau}\log\v{I_\tau} = -2{\rm dist}(D, B).
\label{the formula to get the minimum lengths for CP}
\end{equation}
\end{thm}
\begin{rem}
In Theorem \ref{CP case}, we can relax 
a regularity condition on $\partial{D}$ as we see in condition $(b)$. 
Instead of that, we need convexity for $B$. In the best knowledge of the authors, 
this is a new result. 
\end{rem}

From Theorem \ref{CP case}, we can enclose $D$. 
In what follows, for $p \in \R^3$ and $r > 0$, we put 
$B_{r}(p) = \{x \in \R^3 \vert \v{x - p} < r \}$. 
Take $B$ as $B = B_r(p)$. 
Then $D$ is in the set defined by 
$\v{x - p} > {\rm dist}(D, B_r(p))+r$. 
Hence, if we know $\overline{D} \subset B_{R}(0)$, moving $p$ with $\v{p} > R+r$ so that
$\overline{B_{r}(p)}\cap\overline{D} = \emptyset$, 
we obtain $D \subset \bigcap_{\v{p} > R+r}\{x \in \R^3 \vert \v{x - p} > {\rm dist}(D, B_r(p))+r\}$. 

To obtain Theorem \ref{CP case}, we need 
\begin{align}
I_\tau &\geq \int_{\R^3}(I_3-\gamma(x))\nabla_xv\cdot\nabla_xvdx+O(\tau^{-1}e^{-{\tau}T})
\qquad(\tau \geq 1),
\label{estimate of I_tau for CP -case}
\\
I_\tau &\leq \int_{\R^3}\gamma^{-1/2}(I_3-\gamma(x))\gamma^{-1/2}\nabla_xv\cdot\nabla_xvdx+O(\tau^{-1}e^{-{\tau}T}) \qquad(\tau \geq 1), 
\label{estimate of I_tau for CP +case}
\end{align}
which are given in \cite{IS,  IE, transmission No1}, for example. 
Note that these estimates are shown by only assuming 
$\gamma \in L^\infty(\R^3)$ since $D$ is given by 
$D = \{x \in \R^3 \vert \gamma(x) \neq I_3 \}$.
If $I_\tau \leq 0$, (\ref{estimate of I_tau for CP -case}) yields
$$
\v{I_\tau} = -I_\tau \leq \int_{\R^3}(\gamma(x)-I_3)\nabla_xv\cdot\nabla_xvdx + O(\tau^{-1}e^{-{\tau}T}) \leq C\V{\nabla_xv}_{L^2(D)}^2 + C'\tau^{-1}e^{-{\tau}T}.
$$
If $I_\tau \geq 0$, from (\ref{estimate of I_tau for CP +case}), the same estimate as above
is shown. Thus, we obtain the following key estimates:
\begin{lmm}\label{The key estimate for CP}
(i) There exist constants $C > 0$ and $C' > 0$ such that for $\tau \geq 1$, 
$$
\v{I_\tau} \leq C\V{\nabla_xv}_{L^2(D)}^2 + C'\tau^{-1}e^{-{\tau}T}.
$$
(ii) If $({\rm M})_+$ is assumed, then 
there exist constants $C_0 > 0$ and $C_0' > 0$ such that
\begin{align*}
-I_\tau \geq C_0\V{\nabla_xv}_{L^2(D)}^2 + C_0'\tau^{-1}e^{-{\tau}T}
\qquad(\tau \geq 1).
\end{align*}
(iii) If $({\rm M})_-$ is assumed, then 
there exist constants $C_0 > 0$ and $C_0' > 0$ such that
\begin{align*}
I_\tau \geq C_0\V{\nabla_xv}_{L^2(D)}^2 + C_0'\tau^{-1}e^{-{\tau}T}
\qquad(\tau \geq 1).
\end{align*}
\end{lmm}

Thus, to obtain Theorem \ref{CP case}, it suffices to show the estimate 
of $\V{\nabla_xv}_{L^2(D)}$.
\begin{prp}\label{estimate of nabla_xv for CP}
Suppose that $D$ and $B$ satisfy either $(a)$ or $(b)$ in Theorem \ref{CP case}. 
Then, there exists a constant $C > 0$ such that
$$
C^{-1}\tau^{-7}e^{-2\tau{\rm dist}(D, B)}\leq \V{\nabla_xv(\cdot; \tau)}_{L^2(D)}^2
\leq C\tau^2e^{-2\tau{\rm dist}(D, B)}
\qquad(\tau \geq 1),
$$
where $v$ is the solution of (\ref{Equation of v for CP}).
\end{prp}
From Lemma \ref{The key estimate for CP}, (i) and (ii) of Theorem \ref{CP case} are obvious. 
If we assume $(\text{M})_+$ or $(\text{M})_-$, Lemma \ref{The key estimate for CP} and
Proposition \ref{estimate of nabla_xv for CP} imply
$$
C^{-1}\tau^{-7} -C'\tau^{-1}e^{\tau(-T+2{\rm dist}(D, B))} 
\leq \v{I_\tau}e^{2\tau{\rm dist}(D, B)} 
\leq C\tau^2+C'\tau^{-1}e^{\tau(-T+2{\rm dist}(D, B))}
$$
for some constants $C, C' > 0$ independent of $\tau \geq 1$, which yields 
(\ref{the formula to get the minimum lengths for CP}) in Theorem \ref{CP case}.

\par
 
Proposition \ref{estimate of nabla_xv for CP} shows that the shortest length for this problem is 
$2{\rm dist}(D, B)$. This fact can be understood intuitively. To search $D$, we emit waves from $y \in B$. These waves hit at $x \in D$ and come back to $\xi \in B$. Hence, 
we have to measure until the signals come back. 
In this way, the shortest time appears, which is given by the minimum of 
$\v{y - x}+\v{x - \xi}$ for $x \in \overline{D}$, $y, \xi \in \overline{B}$.

\par
This fact can be found by the representation of 
$\V{\nabla_xv(\cdot; \tau)}_{L^2(D)}^2$. 
In this case, the usual $L^2$-solution $v$ of (\ref{Equation of v for CP}) is given by 
\begin{equation}
v(x; \tau) = \frac{1}{4\pi}\int_{B}\frac{e^{-\tau\v{x - y}}}{\v{x- y}}f(y)dy.
\label{integral kernel representation for CP case}
\end{equation}
Hence, it follows that
\begin{align}
\V{\nabla_xv(\cdot; \tau)}_{L^2(D)}^2 &= \frac{1}{(4\pi)^2}\int_{B}dy\int_{B}d{\xi}f(y)f(\xi)
\int_{D}\frac{e^{-\tau(\v{x - y}+\v{x - \xi})}}{\v{x - y}\v{x - \xi}}
\label{formula of L^2-norm of nabla_x fro CP}
\\&\hskip20mm
\Big(\tau+\frac{1}{\v{x- y}}\Big)\Big(\tau+\frac{1}{\v{x- \xi}}\Big)
\frac{(x - y)\cdot(x - \xi)}{\v{x- y}\v{x- \xi}}dx
\nonumber
\\&
\leq C\V{f}_{L^2(B)}^2\tau^2e^{-2\tau{\rm dist}(D, B)}
\quad(\tau \geq 1),
\nonumber
\end{align}
which gives the upper estimate described in Proposition \ref{estimate of nabla_xv for CP}.

\par

The lower estimate is shown by using similar arguments as in 
\cite{Ikehata-Kawashita2, transmission No1}.
To see the mathematical reason why the shortest length appears, we repeat the argument in \cite{transmission No1} briefly. 
Since we handle some cases with less regularities for the boundaries, 
we begin with estimates for integrals on cones.
\begin{lmm}\label{estimate of integrals on cones}
For any cone $C(a, n, h, \theta)  = \{x \in \R^3 \vert \v{x - a} \leq h, 
(x - a){\cdot}n \geq \v{x - a}\cos\theta \}$, there exists a constant $C > 0$ such that
$$
\int_{C(a, n, h, \theta)}e^{-\tau\v{x-a}}dx \geq C\tau^{-3} \qquad(\tau \geq 1).
$$
\end{lmm}
\noindent
{\quad \it Proof.}
By changing the coordinate, it suffices to take $a = 0$ (i.e. origin) and $n = {}^t(0, 0, 1)$ 
from the beginning. Put $S = \{\omega \in S^2 \vert \v{\omega} = 1, 
\omega_3 \geq \cos\theta \}$ and introduce the polar coordinate $x = r\omega$, 
where $r \geq 0$ and $\omega \in S^2$. By this reduction, 
the cone $C(a, n, h, \theta)$ is given by 
$\{r\omega \vert 0 \leq r \leq h, \omega \in S \}$, which yields
\begin{align*}
\int_{C(a, n, h, \theta)}e^{-\tau\v{x-a}}dx 
&= \int_{0}^{h}\int_{S}e^{-\tau{r}}
r^2 drdS_\omega
= \frac{{\rm Vol}(S)}{\tau^3}\int_{0}^{h\tau}
s^2e^{-s}ds \\
&\geq \frac{{\rm Vol}(S)}{\tau^3}\int_{0}^{h}s^2e^{-s}ds
\end{align*}
for $\tau \geq 1$. This completes the proof of Lemma \ref{estimate of integrals on cones}.
\hfill $\square$
\par
Take $\delta > 0$ sufficiently small so that
\begin{equation}
\frac{(x - y)\cdot(x - \xi)}{\v{x- y}\v{x- \xi}} \geq \frac{1}{2}
\qquad(\v{y - \xi} < 3\delta, y, \xi \in \overline{D}, x \in \overline{B}). 
\label{positivity of amplitude functions}
\end{equation}

\par 

We put $E = \{ (x, y, \xi) \in \overline{D}\times\overline{B}\times\overline{B} \vert \v{y - x}+\v{\xi - x} = 2{\rm dist}(D, B) \}$ and
$E_0 = \{(x, y) \in \overline{D}\times\overline{B} \vert \v{y - x} = {\rm dist}(D, B) \}$. 
Since $B$ and $D$ are bounded, $E_0 \neq \emptyset$, and $(x, y, y) \in E$ for $(x, y) \in E_0$. 
\begin{lmm}\label{characterization of E}
Suppose that $\partial{D}$ has $C^1$ boundary or $B$ is convex. Then $E$ satisfies
$ E = \{ (x, y, \xi) \in \overline{D}\times\overline{B}\times\overline{B} \vert
(x, y) \in E_0, y = \xi \}$. 
\end{lmm}
\noindent
{\quad \it Proof.}
For $(x_0, y_0, \xi_0) \in E$, it follows that
$$
2{\rm dist}(D, B) = \v{y_0 - x_0}+\v{\xi_0 - x_0} \geq \v{y_0 - x_0}+{\rm dist}(D, B),
$$
which yields ${\rm dist}(D, B) \geq \v{y_0 - x_0}$. Hence, we have $(x_0, y_0) \in E_0$. 
We also have $(x_0, \xi_0) \in E_0$ similarly.
\par
We show $y_0 = \xi_0$. First consider the case that $\partial{D}$ is $C^1$. 
Take $(x_0, y_0) \in E_0$ and a $C^1$ curve
$\varphi(t)$ ($\v{t} < \varepsilon$) on $\partial{D}$ with $\varphi(0) = x_0$.
Since $\v{\varphi(t) - y_0}$ takes a minimum at $t = 0$, we have
$$
0 = \frac{d}{dt}\v{\varphi(t) - y_0}\Big\vert_{t = 0}
= \frac{x_0-y_0}{\v{x_0 - y_0}}{\cdot}\varphi'(0), 
$$
which implies that $\frac{y_0-x_0}{\v{y_0 - x_0}} = \nu_{x_0}$, where $\nu_{x_0}$ is the unit 
outer normal vector of $\partial{D}$ at $x_0 \in \partial{D}$. Thus, $y_0$ is given by
$y_0 = \v{y_0 - x_0}\nu_{x_0}+x_0$. Since $(x_0, \xi_0) \in E_0$, we also have
$\xi_0 = \v{\xi_0 - x_0}\nu_{x_0}+x_0$. Noting 
$\v{y_0 - x_0} = \v{\xi_0 - x_0} = {\rm dist}(D, B)$, we obtain $y_0 = \xi_0$. 

\par

Next, assume that $B$ is convex. 
We put $\tilde\varphi(t) = \v{ty_0+(1-t)\xi_0 - x_0}$. 
Since $\overline{B}$ is also convex and $y_0, \xi_0 \in \overline{B}$, 
$ty_0+(1-t)\xi_0 \in \overline{B}$ for $0 \leq t \leq 1$. Hence, it follows that 
\begin{align*}
{\rm dist}(B, D) &\leq \tilde\varphi(t) 
= \v{ty_0+(1-t)\xi_0 - x_0} \leq t\v{y_0 - x_0}+(1-t)\v{\xi_0 - x_0}
\\&
= t{\rm dist}(B, D)+(1-t){\rm dist}(B, D) = {\rm dist}(B, D),
\end{align*}
which yields that $\tilde\varphi(t)={\rm dist}(B, D)$ for $0 \leq t \leq 1$. 
Since 
$${\rm dist}(B, D)^2=t^2\v{y_0-\xi_0}^2+2t(y_0-\xi_0)\cdot(\xi_0-x_0)
+{\rm dist}(B, D)^2 \qquad \text{ ($ 0 \le t \le 1$), } $$
we have $t^2 \v{y_0-\xi_0}^2+2t(y_0-\xi_0)\cdot(\xi_0-x_0)=0$ for any $0 \le t \le 1$.  
This implies that $y_0 = \xi_0$, which completes the proof of Lemma \ref{characterization of E}.
\hfill $\square$
\par
Now we are in a position to show the lower estimate of Proposition
\ref{estimate of nabla_xv for CP}.  
\par
\noindent
{\quad \it Proof.}
Since $\overline{D}$ and $\overline{B}$ are compact, 
for $\delta > 0$ given in (\ref{positivity of amplitude functions}), 
there exist $(x^{(j)}, y^{(j)})$ $\in E_0$ 
($j = 1, \ldots, N$) such that 
$E_0 \subset \bigcup_{j = 1}^{N}B_{\delta}(x^{(j)}){\times}B_{\delta}(y^{(j)})$.
From Lemma \ref{characterization of E}, as a covering of $E$, 
we can take ${\mathcal W} = \bigcup_{j = 1}^{N}
B_{\delta}(x^{(j)}){\times}B_{\delta}(y^{(j)}){\times}B_{\delta}(\xi^{(j)})$. 
Since $E \subset {\mathcal W}$ and 
$\overline{D}{\times}\overline{B}{\times}\overline{B}\setminus{\mathcal W}$ is compact, there exists a constant $c_1 > 0$ such that
\begin{equation}
\v{x - y}+\v{x - \xi} \geq 2{\rm dist}(D, B)+c_1\qquad((x, y, \xi) \in 
\overline{D}{\times}\overline{B}{\times}\overline{B}\setminus{\mathcal W}).
\label{length is longer outside of mathcal W}
\end{equation}
If $(x, y, \xi) \in {\mathcal W}$, $\v{x - x^{(j)}} < \delta $, $\v{y - y^{(j)}} < \delta$ 
and $\v{\xi - y^{(j)}} < \delta$ for some $j \in \{1, 2, \ldots, N\}$, which yields
$\v{y - \xi} \leq \v{y - y^{(j)}}+\v{\xi - y^{(j)}} < 2\delta$. From this and 
(\ref{positivity of amplitude functions}), it follows that 
\begin{equation}
\frac{(x - y)\cdot(x - \xi)}{\v{x- y}\v{x- \xi}} \geq \frac{1}{2}
\qquad((x, y, \xi) \in {\mathcal W}).
\label{positivity estimate near mathcal W}
\end{equation}
Noting (\ref{formula of L^2-norm of nabla_x fro CP}), 
(\ref{length is longer outside of mathcal W}), 
(\ref{positivity estimate near mathcal W}) and (\ref{emission condition}), we obtain
\begin{align}
\V{\nabla_xv(\cdot; \tau)}_{L^2(D)}^2 \geq C\tau^2&\Big\{
\int_{{\mathcal W}\cap(\overline{D}{\times}\overline{B}{\times}\overline{B})}
e^{-\tau(\v{x - y}+\v{x - \xi})}dxdyd{\xi}
\label{estimate of the L^2 norm of the gradient from the below No1}
\\&\hskip10mm
- C'\V{f}_{L^2(B)}^2e^{-2\tau{\rm dist}(D, B)-c_1\tau}\Big\}
\quad(\tau \geq 1), 
\nonumber
\end{align}
where constants $C > 0$ and $C' > 0$ are independent of $\tau$.
Since $(x^{(1)}, y^{(1)})$ attains $\v{x - y} = {\rm dist}(D, B)$, it follows that
$x^{(1)} \in \partial{D}$ and $y^{(1)} \in \partial{B}$. 
From condition (C), there exists a cone $C_{x^{(1)}}$ with vertex $x^{(1)}$ and 
$C_{y^{(1)}}$ with vertex $y^{(1)}$ such that
$C_{x^{(1)}}\setminus\{x^{(1)}\} \subset D$ and 
$C_{y^{(1)}}\setminus\{y^{(1)}\} \subset B$. Note that these cones can be chosen as
$C_{x^{(1)}}\setminus\{x^{(1)}\} \subset D\cap{B_\delta(x^{(1)})}$ and 
$C_{y^{(1)}}\setminus\{y^{(1)}\} \subset B\cap{B_\delta(y^{(1)})}$ if we take the radiuses of
these cones small enough.
Since ${\mathcal W} \supset B_{\delta}(x^{(1)}){\times}B_{\delta}(y^{(1)}){\times}B_{\delta}(y^{(1)})$, from Lemma \ref{estimate of integrals on cones}, it follows that
\begin{align*}
\int_{{\mathcal W}\cap(\overline{D}{\times}\overline{B}{\times}\overline{B})}&
e^{-\tau(\v{x - y}+\v{x - \xi})}dxdyd{\xi}
\geq \int_{C_{x^{(1)}}}dx\left(\int_{C_{y^{(1)}}}e^{-\tau\v{x - y}}dy\right)^2
\\&
\geq \int_{C_{x^{(1)}}}dx\left(\int_{C_{y^{(1)}}}
e^{-\tau\v{x - x^{(1)}}}e^{-\tau\v{x^{(1)}-y^{(1)}}}
e^{-\tau\v{y - y^{(1)}}}dy\right)^2
\\&
= \int_{C_{x^{(1)}}}e^{-2\tau\v{x - x^{(1)}}}dx
\left(\int_{C_{y^{(1)}}}e^{-\tau\v{y - y^{(1)}}}dy\right)^2
e^{-2\tau{\rm dist}(D, B)}
\\&
\geq C\tau^{-9}e^{-2\tau{\rm dist}(D, B)}.
\end{align*}
This estimate and (\ref{estimate of the L^2 norm of the gradient from the below No1}) 
imply the lower estimate of Proposition \ref{estimate of nabla_xv for CP}. 
\hfill $\square$


In view of the above argument, for the time dependent enclosure method,  
it is essential to obtain estimates of the resolvents in order to find the shortest length.
This viewpoint has been already found in the boundary inverse problems for the heat equations by 
using the time dependent enclosure method (cf. 
\cite{Ikehata-Kawashita1, Ikehata-Kawashita2, IK5, Several obstacles for heat IBP}).  
These problems are reduced to 
some elliptic problems with a large parameter, and it is important to show the asymptotic behavior of the kernel of 
the resolvents.  
\par 
We have a long history for the estimates of the resolvents. 
Let $\Omega$ be a bounded domain with $C^\infty$ boundary.  
For the solutions $v(x; \tau)$ of 
$$
\left\{
\begin{array}{lll}
(\tau^2-\triangle)v(x; \tau) = 0 \qquad&\quad\text{in } \Omega,
\\
v(x; \tau) = 1 &\quad\text{on } \partial\Omega,
\end{array}
\right.
$$
to obtain short time asymptotics of the heat kernel, 
Varadhan \cite{Varadhan} shows 
$$
\lim_{\tau \to \infty}\frac{\log\v{v(x; \tau)}}{\tau} = -{\rm dist}(x, \partial\Omega),
$$
where ${\rm dist}(x, \partial\Omega) = \inf_{y \in \partial\Omega}\v{x - y}$. 
In the author's best knowledge, this result gave an original suggestion of the fact
that the shortest length appears from the estimates of the resolvents.
\par
In \cite{IK5, Several obstacles for heat IBP}, 
some estimates of the iterated kernels for the boundary integral equation is needed.  
Nevertheless, usual estimates do not work for controlling 
the errors.  
For seeking the shortest length, 
a precise estimate of those kernels 
is developed in \cite{Kernel estimates}.  
\par
Using the precise estimates of the integral kernels, we can also handle the case 
where there is no signal in some part of the boundary. 
Take a bounded set $D$ with $\overline{D} \subset \Omega$. 
Consider the following reduced problems:
$$
\left\{
\begin{array}{lll}
(\tau^2-\triangle)v(x; \tau) = 0 \qquad&\quad\text{in } \Omega\setminus\overline{D},
\\
\partial_{\nu}v(x; \tau) = 1 &\quad\text{on } \partial\Omega, 
\\
\partial_{\nu}v(x; \tau) = 0 &\quad\text{on } \partial{D}, 
\end{array}
\right.
$$
where $\partial_{\nu}$ is the normal derivative of the boundary.
Since there is no signal on $\partial{D}$, 
only the reflected signals 
of original source given in $\partial\Omega$ come back to the inside. 
Hence, in this case, 
${\rm dist}(x, \partial\Omega)$ is the shortest length from the boundary to
$x \in \Omega\setminus{D}$, which is justified in \cite{Varadhan Extended}
(cf. Theorem 1.3 in \cite{Varadhan Extended}).
To obtain this result, we essentially use the precise estimate given 
in \cite{Kernel estimates}. 
Thus, these problems described above are linked to each other by finding 
the shortest lengths corresponding to the problems.

\setcounter{equation}{0}
\section{Inverse problems in a two-layered medium} 
\label{Inverse problems in two-layered medium}

In this section, we consider inverse problems for detecting inclusions embedded in the below of a two-layered medium whose layer is given by a flat transmission boundary (cf. \cite{transmission No1, transmission No2}). 
We put $\R^3_{\pm} = \{ x = (x_1,x_2,x_3) \in \R^3 \,\vert\, \pm x_3 > 0 \}$. Assume that
the propagation speed of the wave in $\R^3_{\pm}$
is given by $\sqrt{\gamma_\pm}$, where $\gamma_\pm > 0$ are constants 
with $\gamma_+ \neq \gamma_-$. 
We call $\R^3_+$ 
(resp. $\R^3_{-}$) the upper (resp. lower) side of the flat transmission boundary 
$\partial\R^3_{\pm}$. Note that the waves propagating this two-layered medium 
is governed by 
\begin{equation}
\left\{
\begin{array}{ll}
\displaystyle
(\partial_t^2-L_{\gamma_0}) u_0 = 0 \quad&  \qquad\text{in}\, (0, T)\times\R^3,
\\
\displaystyle
u_0(0, x) = 0,  \quad 
\partial_tu_0(0, x) = f(x) & \qquad\text{on}\,\R^3,
\end{array}
\right.
\label{equation of free case}
\end{equation}
where $\gamma_0(x) = \gamma_{\pm}I_3$  for $\pm{x_3} > 0$ and
$L_{\gamma_0}u = {\rm div}(\gamma_0{\nabla}u)$. 

\par

We assume that $\overline{D} \subset \R^3_-$ 
and $\R^3_-\setminus{D}$ is connected. 
Thus, inclusions $D$ are in the lower side. 
We put symmetric matrix valued $\gamma = (\gamma_{ij}(x)) \in L^\infty(\R^3)$ satisfying
$\gamma(x) = \gamma_+I_3$ for $x_3 > 0$, $\gamma(x) = \gamma_-I_3$ in $\R^3_-\setminus{D}$ and either $(\text{M})_{\gamma_0,+}$ or $(\text{M})_{\gamma_0,-}$, where
$(\text{M})_{\gamma_0,+}$ : $\gamma$ satisfies $\inf_{x \in D}\gamma(x) > \gamma_- I_3$, $(\text{M})_{\gamma_0,-}$ : $\gamma$ satisfies $\sup_{x \in D}\gamma(x) < \gamma_- I_3$.
The waves propagating this medium with inclusions satisfy the same equations as 
(\ref{equation with inclusion}) with $\gamma$ given in this section.

\par

In the two-layer case, it is a natural restriction that signals can not be emitted in the lower side. Hence, we take an open set $B$ with $\overline{B} \subset \R^3_+$ and 
the data $f$  so that they satisfy the emission condition (\ref{emission condition}) 
in the same equations as (\ref{equation with inclusion}) with $\gamma$.

\par

We introduce the indicator function $I_\tau$ by (\ref{Indicator function for CP}), where
$w$ is given by (\ref{def of w}) for the solution $u$ of the two-layered medium with inclusions.
In this case, $v$ should be taken as the 
$L^2$-solution of 
\begin{equation}
(L_{\gamma_0}-\tau^2)v(x; \tau)+f(x) = 0 \qquad\text{in } \R^3
\label{Equation of v for transmission case}
\end{equation}
since the medium without inclusions is the two-layered medium with the flat transmission 
boundary.

\par

Similarly to Theorem \ref{CP case}, we can find the shortest length in this case, which is
given by
\begin{equation*}
\displaystyle
l(D, B)=\inf_{x \in D,\, y \in B}l(x, y),
\end{equation*}
where
\begin{align}
l(x, y) &= \inf_{z' \in \R^2}l_{x, y}(z') \qquad(x \in \R^3_-, y \in \R^3_+), 
\nonumber
\\
l_{x, y}(z') &= \frac{1}{\sqrt{\gamma_-}}\v{\tilde{z}'-x}
+\frac{1}{\sqrt{\gamma_+}}\v{\tilde{z}'-y}
\quad(\tilde{z}' = (z_1, z_2, 0), z' = (z_1, z_2)).
\label{the path length for Snell's law}
\end{align}
Note that $l(x, y)$ (resp. $l(D, B)$) is called the {\it optical distance} between $x$ and $y$ 
(resp. the inclusion $D$ and the set $B$).

\par

The optical distance between $x \in \R^3_-$ and $y \in \R^3_+$ is determined by Snell's law. 
As is in Lemma 4.1 of \cite{transmission No1}, for arbitrary $x$ and $y \in \R^3$ 
with $x_3 < 0$ and $y_3 > 0$, there exists a unique point $z'(x, y) \in \R^2$ 
satisfying $l(x, y) = l_{x, y}(z'(x, y))$, and the point $z'(x, y)$ is on the
line segment $x'y'$ and $C^\infty$ for $x$ and $y \in \R^3$ with 
$x_3 < 0$ and $y_3 > 0$. 
Since $z' = z'(x, y)$ is a unique point attaining $l(x, y) = \inf_{z' \in \R^2}l_{x, y}(z')$,
$\tilde{z}'(x, y) = (z'(x, y), 0)$ satisfies $\partial_{z'}l_{x, y}(z') = 0$, which yields
\begin{align}
\frac{1}{\sqrt{\gamma_-}}\frac{z'(x, y) - x'}{\v{\tilde{z}'(x, y)-x}}
+\frac{1}{\sqrt{\gamma_+}}\frac{z'(x, y) - y'}{\v{\tilde{z}'(x, y)-y}} = 0. 
\label{Snell's law No1}
\end{align}
We define $0 \leq \theta_\pm < \pi/2$ by
\begin{align}
\sin\theta_- = \frac{\v{{z}'(x, y)- x'}}{\v{\tilde{z}'(x, y)-x}}, \qquad
\sin\theta_+ = \frac{\v{{z}'(x, y)- y'}}{\v{\tilde{z}'(x, y)-y}}. 
\label{how to define sin theta_pm}
\end{align}
Then, the relation (\ref{Snell's law No1}) implies
$$
\frac{\sin\theta_-}{\sqrt{\gamma_-}}\frac{z'(x, y) - x'}{\v{z'(x, y) - x'}}
+ \frac{\sin\theta_+}{\sqrt{\gamma_+}}\frac{z'(x, y) - y'}{\v{z'(x, y) - y'}}
= 0.
$$
This means that $z'(x, y) \in \R^2$ is on the line segment $x'y'$ on $\R^2$, and 
\begin{equation}
\frac{\sin\theta_-}{\sqrt{\gamma_-}} = \frac{\sin\theta_+}{\sqrt{\gamma_+}},
\label{the Snell's law}
\end{equation}
which is just Snell's law.

\par 

The results of two-layer case are formally same as Theorem \ref{CP case}, 
which is described by replacing ${\rm dist}(D, B)$ with $ l(D, B)$: the shortest length 
in this case. 

\begin{thm}\label{transmission case}
Suppose that $D$ and $B$ introduced in this section 
satisfy either $(a)$ or $(b)$:\\
$(a)$ $D$ has $C^1$ boundary $\partial{D}$ and $B$ satisfies condition $(C)$, \\
$(b)$ $D$ and $B$ satisfy condition $(C)$ and $B$ is convex. \\
Then, the indicator function (\ref{Indicator function for CP}) defined by the solution
$v$ of (\ref{Equation of v for transmission case}) satisfies 
\par\noindent
(i) for $T < 2l(D, B)$, $\lim_{\tau \to \infty}e^{{\tau}T}I_\tau = 0$, 
\par\noindent
(ii) for $T > 2l(D, B)$, $\lim_{\tau \to \infty}e^{{\tau}T}I_\tau = \mp\infty$ 
if $(\text{M})_{\gamma_0, \pm}$ is assumed, respectively.
\par
Further, suppose $T > 2l(D, B)$, and either $(\text{M})_{\gamma_0, +}$ or $(\text{M})_{\gamma_0, -}$. Then
\begin{equation}
\lim_{\tau \to \infty}\frac{1}{\tau}\log\v{I_\tau} = -2l(D, B).
\nonumber
\end{equation}
\end{thm}


Take $p \in \R^3_+$ and $r > 0$ so that $\overline{B_r(p)} \subset \R^3_+$. Then, 
Theorem \ref{transmission case} implies $D$ is in the set $\{ x \in \R^3_- \vert 
l(x, p) > l(D, B_r(p))+r/\sqrt{\gamma_+}\}$. Since 
$l(x, p) \leq \v{x - \tilde{z}'}/\sqrt{\gamma_-} + \v{p - \tilde{z}'}/\sqrt{\gamma_+}$ for
any $z' \in \R^2$, it follows that
$$
D \subset\bigcap_{p \in \R^3_+, r < p_3}\bigcap_{z' \in \R^2}
\Big\{x \in \R^3_- \vert  \v{x - \tilde{z}'} > \sqrt{\gamma_-}\Big(l(D, B_r(p))
+\frac{r-\v{p - \tilde{z}'}}{\sqrt{\gamma_+}}\Big)\Big\}.
$$ 
Thus, in this case, $D$ can be enclosed from the upper side.
This result is reasonable with this setting, because 
we only can emanate signals to the inclusion $D$ and 
catch the reflected waves.

\par

To show Theorem \ref{transmission case}, we basically follow the arguments 
in section \ref{no layered case}. Similarly to the case of the Cauchy problems in 
section \ref{no layered case}, Theorem \ref{transmission case} is given by the following estimates:
\begin{prp}\label{estimate of nabla_xv for transmission case}
Suppose that $D$ and $B$ introduced in this section satisfy either 
$(a)$ or $(b)$ in Theorem \ref{transmission case}.  
Then, there exists a constant $C > 0$ such that
$$
C^{-1}\tau^{-7}e^{-2{\tau}l(D, B)}\leq \V{\nabla_xv(\cdot; \tau)}_{L^2(D)}^2
\leq C\tau^2e^{-2{\tau}l(D, B)}
\qquad(\tau \geq 1),
$$
where $v$ is the solution of (\ref{Equation of v for transmission case}).
\end{prp}

In this case, the solution $v$ of (\ref{Equation of v for transmission case}) is given by
$$
v(x; \tau) = \int_{\R^3}\Phi_\tau(x, y)f(y)dy,
$$
where $\Phi_\tau(x, y)$ is the fundamental solution of (\ref{Equation of v for transmission case}). 
Hence, it follows that
\begin{align}
\int_{D}\v{\nabla_xv(x; \tau)}^2dx
= \int_{B}dy\int_{B}d{\xi}f(y){f(\xi)}
\int_{D}\nabla_x\Phi_\tau(x, y)\cdot{\nabla_x\Phi_\tau(x, \xi)}dx. 
\label{L^2 norms of nabla_xv in D}
\end{align}
Thus, the aim of this problem is to obtain asymptotic behavior of $\nabla_x\Phi_\tau(x, y)$ as 
$\tau \to \infty$.
Take open sets $B'$ and $D'$ with $\overline{B'} \subset \R^3_+$, 
$\overline{D'} \subset \R^3_-$, 
$\overline{B} \subset B'$ and $\overline{D} \subset D'$.
\begin{prp}\label{Asymptotic behaviour of Phi_tau(x, y)}
Assume that $\gamma_+ \neq \gamma_-$. Then for $k = 0, 1$, we have 
\begin{align}
\nabla_x^k\Phi_\tau(x, y) = 
\frac{{e}^{-{\tau}l(x, y)}}{8\pi\gamma_+\gamma_-
\sqrt{{\rm det}H(x, y)}}
\Big(\frac{-\tau}{\sqrt{\gamma_-}}\Big)^k
\Big(\sum_{j = 0}^{N}{\tau^{-j}}\Phi_{j}^{(k)}(x, y) 
+ Q_{N, \tau}^{(k)}(x, y)\Big), 
\nonumber
\end{align}
where $H(x, y) = {\rm Hess}(l_{x, y})(z'(x, y))$ is the Hessian of $l_{x, y}$ given by
(\ref{the path length for Snell's law}) at $z' = z'(x, y)$, 
$\Phi_{j}^{(k)}(x, y)$ $(k = 0, 1)$ are 
$C^\infty$ in $\overline{D'}\times\overline{B'}$,
for any $N \in {\N}\cup\{0\}$, $Q_{N, \tau}^{(k)}(x, y)$ $(k = 0, 1)$
are continuous in $\overline{D'}\times\overline{B'}$
with a constant $C_N > 0$ satisfying
$$
\v{Q_{N, \tau}^{(0)}(x, y)} + \v{Q_{N, \tau}^{(1)}(x, y)} \leq C_N\tau^{-(N+1)} \qquad(x \in \overline{D'}, y \in \overline{B'}, \tau \geq 1).
$$
Moreover, $\Phi_{0}^{(k)}(x, y) $ $(k = 0, 1)$ are given by
\begin{align*}
\Phi_{0}^{(0)}(x, y)
&= \frac{E_{0}(x - \tilde{z}'(x, y))}{\v{x - \tilde{z}'(x, y)}\v{\tilde{z}'(x, y)-y}},
\intertext{and}
\Phi_{0}^{(1)}(x, y) 
&= \Phi_{0}^{(0)}(x, y)\frac{x - \tilde{z}'(x, y)}{\v{x - \tilde{z}'(x, y)}},
\end{align*}
where
\begin{align}
E_{0}(x - \tilde{z}') 
&= \frac{4\sqrt{\gamma_-}\v{x_3}\sqrt{a_0^2\v{x - \tilde{z}'}^2- \v{x' - z'}^2}}
{\v{x - \tilde{z}'}\big(\sqrt{a_0^2\v{x - \tilde{z}'}^2- \v{x' - z'}^2}+a_0^2\v{x_3}\big)}, 
\quad a_0=\sqrt{\frac{\gamma_-}{\gamma_+}}.
\nonumber
\end{align}
\end{prp}

The most important part for proving Theorem \ref{transmission case} is to show
Proposition \ref{Asymptotic behaviour of Phi_tau(x, y)}, which is obtain 
in \cite{transmission No1} for $\gamma_+ < \gamma_-$ and 
in \cite{transmission No2} for $\gamma_+ > \gamma_-$. 
The difference between the cases $\gamma_+ < \gamma_-$ and $\gamma_+ > \gamma_-$ is caused by the total reflection phenomena. 
Signals from $B$ can reach $D$ as waves refracted by the transmission boundary. 
These waves are reflected by the inclusion $D$, 
and come back to the transmission boundary. 
The problems are whether the refracted waves of these waves 
arriving at the transmission boundary 
go back to $B$ or not. 

\par

If $\gamma_+ < \gamma_-$, $\theta_\pm$ defined by (\ref{how to define sin theta_pm})
satisfy $\theta_- > \theta_+$, which yields that every wave coming from the lower side goes to the upper side. 
This means that the waves from the lower side do not cause the total reflection. 
In this case, if a wave reflects at $x \in D$, passes $\tilde{z}' \in \partial\R^3_+ $ and
reaches $y \in B$, it takes time $l_{x, y}(z')$ given by (\ref{the path length for Snell's law}). 
Hence, $l(x, y) = \inf_{z' \in \R^3_+}l_{x, y}(z')$ gives the time it takes waves to travel 
from $x$ to $y$. 
This fact affects on the proof of Proposition \ref{Asymptotic behaviour of Phi_tau(x, y)}, 
given in \cite{transmission No1}.
We can have a good complex integral representation of $\Phi_\tau(x, y)$ 
for $x_3 < 0$.

\par

On the other hand, if $\gamma_+ > \gamma_-$, we have $\theta_- < \theta_+$ and
the total reflection phenomena make the problems more complicated. 
We put $\theta_0 \in (0, \pi/2)$ defined by $\sin\theta_0 = a_0=\sqrt{\gamma_-/\gamma_+}$,
which is called the critical angle. 
If $\theta_- < \theta_0$, $\theta_+ \in [0, \pi/2)$ is 
determined by Snell's law (\ref{the Snell's law}). 
This means that the waves reflected by 
the inclusion $D$ are refracted, and propagate into the upper side. 
If $\theta_- \geq \theta_0$, (\ref{the Snell's law}) does not make sense except the case
$\theta_- = \theta_0$. 
In this case, waves propagating along the transmission boundary $\partial\R^3_-$ appear. 
These waves are called \lq\lq evanescent waves'', which propagate along $\partial\R^3_+$ 
and emit the body waves into the upper side. 
Corresponding to this fact, $l_{x, y}(z')$ is different 
from the time it takes waves to travel from $x \in \R^3_-$ 
to $y \in \R^3_+$. As is in \cite{transmission No2}, 
this time is estimated by 
\begin{align}
\tilde{l}_{x, y}(z') = 
\begin{cases} 
l_{x, y}(z') & (z' \in {\mathcal U}_{1}(x)), \\[2mm]
\displaystyle
\frac{\v{x_3}\cos\theta_0}{\sqrt{\gamma_-}}
+\frac{\v{x' - z'}+\v{\tilde{z}' - y}}{\sqrt{\gamma_+}} & 
(z' \in \R^2{\setminus}{\mathcal U}_{1}(x)),
\end{cases}
\nonumber
\end{align}
where $ {\mathcal U}_{\delta}(x) = \{\, z' \in \R^2 \,\vert\, \v{x' - z'} 
< a_0{\delta}\v{x - \tilde{z}'} \,\}$
for $0 < \delta < a_0^{-1}$.  
Note that for $z' \in {\mathcal U}_{1}(x)$, 
$\tilde{l}_{x, y}(z') = l_{x, y}(z')$ is the exact time.  
When $z' \in {\mathcal U}_{1}(x)$, we have $\theta_- < \theta_0$ and 
$\tilde{l}_{x, y}(z') = l_{x, y}(z')$ for the trajectories corresponding to the usual 
refracted waves.
Note that the function $\tilde{l}_{x, y}(z')$ is determined in the steps for getting
asymptotic behavior of $\Phi_\tau(x, y)$. 
When $z' \in \R^2{\setminus}{\mathcal U}_{1}(x)$, we have $\theta_- \ge \theta_0$ and 
$\Phi_\tau(x, y)$ has a different complex integral representation
from that for the case $\gamma_+ < \gamma_-$. This difference leads the function 
$\tilde{l}_{x, y}(z')$ (cf. sections 2 and 3 of \cite{transmission No2}).

\par 

To explain the meaning of $\tilde{l}_{x, y}(z')$ for $z' \in \R^2{\setminus}{\mathcal U}_{1}(x)$, 
we take $z_0' \in \R^2$ on the line segment $x'z'$ with 
$ \v{x' - z_0'}/\v{x - \tilde{z}_0'} = \sin\theta_0$ and
$\v{x' - z'} = \v{x' - z_0'}+\v{z_0' - z'}$. Note that this $z_0'$ can be chosen
since $\v{x' - z'}/\v{x - \tilde{z}'} > \sin\theta_0$ for 
$z' \in \R^2{\setminus}{\mathcal U}_{1}(x)$. Hence, for any $y \in \R^3_+$, 
it follows that
\begin{align}
\tilde{l}_{x, y}(z')
&=\frac{\cos\theta_0 }{\sqrt{\gamma_-}}
\frac{\v{x_3}}{\v{x-\tilde{z}'_0}}
\v{x-\tilde{z}'_0}
+\frac{\v{x'-z_0'}+\v{z_0'-z'}}{\sqrt{\gamma_+}} 
+\frac{\v{\tilde{z}'-y}}{\sqrt{\gamma_+}}
\nonumber \\
&=\frac{\cos^2\theta_0 }{\sqrt{\gamma_-}}
\v{x-\tilde{z}'_0}
+\frac{\sqrt{\gamma_-}}{\sqrt{\gamma_+}}
\frac{\v{x'-z_0'}}{\v{x-\tilde{z}'_0}}
\frac{\v{x-\tilde{z}'_0}}{\sqrt{\gamma_-}}
+\frac{\v{z_0'-z'}}{\sqrt{\gamma_+}}
+\frac{\v{\tilde{z}'-y}}{\sqrt{\gamma_+}}
\nonumber\\
&=\frac{\v{x-\tilde{z}'_0}}{\sqrt{\gamma_-}}
+\frac{\v{z_0'-z'}+\v{\tilde{z}'-y}}{\sqrt{\gamma_+}}. 
\nonumber
\end{align}
From this result, we can say that for 
$z' \in \R^2{\setminus}{\mathcal U}_{1}(x)$, 
the optical distance between $x$ and $\tilde{z}'$ is not
$\v{x - \tilde{z}'}/\sqrt{\gamma_-}$, and is given by
$\v{x-\tilde{z}'_0}/\sqrt{\gamma_-}+\v{z_0'- z'}/\sqrt{\gamma_+}$.

\par

For $z' \in \R^2{\setminus}{\mathcal U}_{1}(x)$, it follows that
$$
\v{x - \tilde{z}'}^2 = (\v{x-\tilde{z}'_0} + a_0\v{z_0'- z'})^2
+ \frac{x_3^2\v{z' - z_0'}^2}{\v{x - \tilde{z}_0'}^2}, 
$$
which yields $ \v{x-\tilde{z}'_0} + a_0\v{z_0'- z'} < \v{x - \tilde{z}'}$.
Hence, for $z' \in \R^2{\setminus}{\mathcal U}_{1}(x)$,  
$\tilde{l}_{x, y}(z')$ is shorter than $l_{x, y}(z')$ since 
\begin{align*}
\frac{\v{x-\tilde{z}'_0}}{\sqrt{\gamma_-}}+\frac{\v{z_0'- z'}}{\sqrt{\gamma_+}}
= \frac{1}{\sqrt{\gamma_-}}(\v{x-\tilde{z}'_0} + a_0\v{z_0'- z'})
< \frac{\v{x- \tilde{z}'}}{\sqrt{\gamma_-}}.
\end{align*}
This means that if $\gamma_+ > \gamma_-$, the shortest time from $x \in \R^3_-$ to
$y \in \R^3_+$ should be given by $\inf_{z' \in \R^2}\tilde{l}_{x, y}(z')$. 
\begin{lmm}\label{the shortest length for gamma_+ > gamma_-}
In the case of $\gamma_+ > \gamma_-$, for any $x \in \R^3_-$ and $y \in \R^3_+$, we have 
$\inf_{z' \in \R^2}\tilde{l}_{x, y}(z') = l(x, y)$. The minimum is attained
at only a unique point $z' = z'(x, y)$, which is the same point as the minimizer of 
$\inf_{z' \in \R^2}l_{x, y}(z')$.
\end{lmm} 
For a proof, see Lemma 3.1 of \cite{transmission No2}. From this fact, as in 
Proposition \ref{Asymptotic behaviour of Phi_tau(x, y)}, in the both 
cases $\gamma_+ < \gamma_-$ and $\gamma_+ > \gamma_-$, the asymptotic behavior of 
$\Phi_\tau(x, y)$ are the same form though the reasons are different from each other. 
This is the outline for the two-layer case developed by 
\cite{transmission No1, transmission No2}.

\setcounter{equation}{0}
\section{Estimates of $\V{\nabla_xv(\cdot; \tau)}_{L^2(D)}^2$
for less regularities on $\partial{D}$}
\label{how to show the case of less regularities}

In this section, we show Proposition \ref{estimate of nabla_xv for transmission case}.
From Proposition \ref{Asymptotic behaviour of Phi_tau(x, y)}, 
the estimate from the above in
Proposition \ref{estimate of nabla_xv for transmission case} 
is obtained in the same way as 
Proposition \ref{estimate of nabla_xv for CP}. 
Thus, the problem is to obtain the estimates from the below.

\par 

In \cite{transmission No1, transmission No2}, 
Proposition \ref{estimate of nabla_xv for transmission case} is shown 
by assuming that $\partial{D}$ is $C^1$.
In this case, if $x_0 \in \overline{D}$ and 
$y_0 \in \overline{B}$ attain $l(D, B) = \inf_{x \in \overline{D}, y \in \overline{B}}l(x, y)$,
then, $x_0 \in \partial{D}$ and $y_0 \in \partial{B}$, and $x_0$ satisfies
$$
\nu_{x_0} = \frac{\tilde{z}'(x_0, y_0) - x_0}{\v{\tilde{z}'(x_0, y_0) - x_0}}, \qquad
$$
where $\nu_{x_0}$ is the unit outer normal of $D$ at $x_0 \in \partial{D}$.
As in Lemma 5.1 of \cite{transmission No1} 
this fact is used to show $y_0 = \xi_0$ if $l(x_0, y_0) = l(x_0, \xi_0) = l(D, B)$,  
which corresponds to Lemma \ref{characterization of E}.

\par 

Similarly to Theorem \ref{CP case}, we can relax the assumption for regularities 
of $\partial{D}$. 
We begin with showing the following lemma:
\begin{lmm}\label{uniqueness of points attaining l(D, B)}
Suppose that $\partial{D}$ is $C^1$ or $B$ is convex. Then, for any 
$x_0 \in \overline{D}$, $y_0 \in \overline{B}$ and $\xi_0 \in \overline{B}$ satisfying
$l(x_0, y_0) = l(x_0, \xi_0) = l(D, B)$, we have $y_0 = \xi_0$.
\end{lmm}
\noindent
{\quad \it Proof.}
As described above, if $\partial{D}$ is $C^1$, 
Lemma \ref{uniqueness of points attaining l(D, B)} is shown as
Lemma 5.1 of \cite{transmission No1}. Here, we consider the case
that $B$ is convex.  
We put 
\begin{align*}
\tilde\varphi(t) = &\frac{\v{t\tilde{z}'(x_0, y_0)+(1-t)\tilde{z}'(x_0, \xi_0)- x_0}}{\sqrt{\gamma_-}} \\
&+\frac{\v{t\tilde{z}'(x_0, y_0)+(1-t)\tilde{z}'(x_0, \xi_0)- (ty_0+(1-t)\xi_0)}}{\sqrt{\gamma_+}}.
\end{align*}
Since $\overline{B}$ is also convex and $y_0, \xi_0 \in \overline{B}$, 
$ty_0+(1-t)\xi_0 \in \overline{B}$ for $0 \leq t \leq 1$. 
Hence, from Lemma \ref{the shortest length for gamma_+ > gamma_-} it follows that 
\begin{align*}
l(D,B) \leq \tilde\varphi(t) 
=&\frac{1}{\sqrt{\gamma_-}}\v{t(\tilde{z}'(x_0, y_0)-x_0)+(1-t)(\tilde{z}'(x_0, \xi_0)- x_0)} \\
&+\frac{1}{\sqrt{\gamma_+}}\v{t(\tilde{z}'(x_0, y_0)-y_0)+(1-t)(\tilde{z}'(x_0, \xi_0)-\xi_0)} \\
\leq &t\left\{ 
\frac{1}{\sqrt{\gamma_-}}\v{\tilde{z}'(x_0, y_0)-x_0}
+\frac{1}{\sqrt{\gamma_+}}\v{\tilde{z}'(x_0, y_0)-y_0} \right\} \\
&+(1-t)\left\{\frac{1}{\sqrt{\gamma_-}}\v{\tilde{z}'(x_0, \xi_0)- x_0}
+\frac{1}{\sqrt{\gamma_+}}\v{\tilde{z}'(x_0, \xi_0)-\xi_0} \right\} \\
\leq & t \, l(D,B)+(1-t)\, l(D,B) =l(D,B),
\end{align*}
which yields that $\tilde\varphi(t)$ is constant for $0 \leq t \leq 1$. 
Thus, we have
\begin{align*}
0 = &\frac{d}{dt}\tilde\varphi(t) \Bigm\vert_{t \to +0} \Bigm. \\
= &\frac{1}{\sqrt{\gamma_-}}\frac{(\tilde{z}'(x_0, \xi_0)-x_0)
\cdot(\tilde{z}'(x_0, y_0)- \tilde{z}'(x_0, \xi_0))}
{\v{\tilde{z}'(x_0, \xi_0)-x_0}}\\
&+\frac{1}{\sqrt{\gamma_+}}\frac{(\tilde{z}'(x_0, \xi_0)-\xi_0)
\cdot(\tilde{z}'(x_0, y_0)- \tilde{z}'(x_0, \xi_0)-y_0+\xi_0)}
{\v{\tilde{z}'(x_0, \xi_0)-\xi_0}}\\
= &\frac{1}{\sqrt{\gamma_-}}\Big\{\v{\tilde{z}'(x_0, y_0)-x_0}p(\xi_0){\cdot}p(y_0)
-\v{\tilde{z}'(x_0, \xi_0)-x_0}\Big\}
\\
&
+ \frac{1}{\sqrt{\gamma_+}}\Big\{\v{\tilde{z}'(x_0, y_0)-y_0}q(\xi_0){\cdot}q(y_0)
-\v{\tilde{z}'(x_0, \xi_0)-\xi_0}\Big\}, 
\end{align*}
where $p(y) = \frac{\tilde{z}'(x_0, y)-x_0}{\v{\tilde{z}'(x_0, y)-x_0}}$ and
$q(y) = \frac{\tilde{z}'(x_0, y)- y}{\v{\tilde{z}'(x_0, y)- y}}$.
Since 
$l(D, B) 
= \frac{1}{\sqrt{\gamma_-}}\v{\tilde{z}'(x_0, \xi_0)-x_0} 
+\frac{1}{\sqrt{\gamma_+}}\v{\tilde{z}'(x_0, \xi_0)-\xi_0}
=  \frac{1}{\sqrt{\gamma_-}}\v{\tilde{z}'(x_0, y_0)-x_0} 
+\frac{1}{\sqrt{\gamma_+}}\v{\tilde{z}'(x_0, y_0)-y_0}
$,
we have
$$
\frac{\v{\tilde{z}'(x_0, y_0)-x_0}}{\sqrt{\gamma_-}}\Big\{p(\xi_0){\cdot}p(y_0)-1\Big\}
+ 
\frac{\v{\tilde{z}'(x_0, y_0)-y_0}}{\sqrt{\gamma_+}}\Big\{q(\xi_0){\cdot}q(y_0)-1\Big\}
= 0.
$$
Using $p(\xi_0){\cdot}p(y_0) \leq 1$ and $q(\xi_0){\cdot}q(y_0) \leq 1$, we obtain
$p(\xi_0){\cdot}p(y_0) = q(\xi_0){\cdot}q(y_0) = 1$, 
which yields $p(\xi_0) = p(y_0)$ since
$$
\v{p(\xi_0) - p(y_0)}^2 
= \v{p(y_0)}^2+\v{p(\xi_0)}^2-2p(\xi_0){\cdot}p(y_0) = 0. 
$$
Hence, it follows that
$$
\frac{z'(x_0, y_0)-x_0'}{\v{\tilde{z}'(x_0, y_0)-x_0}}
= \frac{z'(x_0, \xi_0)-x_0'}{\v{\tilde{z}'(x_0, \xi_0)-x_0}}, \qquad
\frac{(x_0)_3}{\v{\tilde{z}'(x_0, y_0)-x_0}}
= \frac{(x_0)_3}{\v{\tilde{z}'(x_0, \xi_0)-x_0}},
$$
which yields $\tilde{z}'(x_0, y_0) = \tilde{z}'(x_0, \xi_0)$
since $(x_0)_3 < 0$.
Put $\psi(t)=\v{ty_0+(1-t)\xi_0 - \tilde{z}'(x_0, y_0)}/\sqrt{\gamma_+}$. 
From the facts $z'(x_0, y_0) = z'(x_0, \xi_0)$ and $\tilde\varphi(t)=l(D,B)$, 
it follows that 
$$
\psi(t)= l(D, B) - \v{x_0 - \tilde{z}'(x_0, y_0)}/\sqrt{\gamma_-}, 
$$
which implies that $\psi(t)$ is independent of $t$. 
Hence, $y_0 = \xi_0$ is shown in a similar way to Lemma \ref{characterization of E}.
\hfill $\square$

\par

Once we obtain Lemma \ref{uniqueness of points attaining l(D, B)}, 
the proof of Proposition \ref{estimate of nabla_xv for transmission case} can be 
obtained in the similar way to the case of Theorem \ref{CP case}, 
and the argument is developed in section 5 of \cite{transmission No1}. 
Hence, we only give the outline of the proof.

{\quad \it Proof.}
From Proposition \ref{Asymptotic behaviour of Phi_tau(x, y)}, it follows that
\begin{align*}
\nabla_x\Phi_\tau(x, y)\cdot{\nabla_x\Phi_\tau(x, \xi)}
&= \frac{\tau^2e^{-\tau(l(x, y)+l(x, \xi))}}{64\pi^2\gamma_+^2\gamma_-^3\sqrt{{\rm det}H(x, y){\rm det}H(x, \xi)}}
\\&\hskip10mm
\Big(h_0(x, y, \xi)\Phi_{0}^{(0)}(x, y)\Phi_{0}^{(0)}(x, \xi)
+\tau^{-1}h_1(x, y, \xi, \tau)\Big),
\end{align*}
where $h_0(x, y, \xi)$ is given by
$$
h_0(x, y, \xi) = \frac{x - \tilde{z}'(x, y)}{\v{x - \tilde{z}'(x, y)}}
\cdot\frac{x - \tilde{z}'(x, \xi)}{\v{x - \tilde{z}'(x, \xi)}}
$$
and $h_1(x, y, \xi, \tau)$ is a continuous function satisfying 
$$
\v{h_1(x, y, \xi, \tau)} \leq C \qquad(x \in \overline{D}, y, \xi \in \overline{B}, \tau \geq 1) 
$$
for some fixed constant.

\par

Since $h_0(x, y, y) = 1$, there exists a constant $\delta > 0$ such that 
$$
h_0(x, y, \xi) \geq 1/2 \qquad(\v{y - \xi} < 3\delta, x \in \overline{D}, y, \xi \in \overline{B}),
$$
and $\overline{B_{4\delta}(y)} \subset \R^3_+$ ($y \in \overline{B}$). 
From Lemma \ref{uniqueness of points attaining l(D, B)}, and compactness of 
$\overline{D}$ and $\overline{B}$, 
for this $\delta$, there exist $(x_0^{(j)}, y_0^{(j)}) \in E_0$ such that 
$E \subset \bigcup_{j = 1}^N
B_{\delta}(x_0^{(j)}){\times}B_{\delta}(y_0^{(j)}){\times}B_{\delta}(y_0^{(j)})$.
We put ${\mathcal W} = \bigcup_{j = 1}^N
B_{\delta}(x_0^{(j)}){\times}B_{\delta}(y_0^{(j)}){\times}B_{\delta}(y_0^{(j)})$.
The argument showing (\ref{length is longer outside of mathcal W}) 
and (\ref{positivity estimate near mathcal W}) implies that 
there exists a constant $c_1 > 0$ such that
\begin{align*}
l(x, y)+l(x, \xi) &\geq 2l(D, B)+c_1\qquad((x, y, \xi) \in 
\overline{D}{\times}\overline{B}{\times}\overline{B}\setminus{\mathcal W}),
\\
h_0(x, y, \xi) &\geq 1/2 \qquad((x, y, \xi) \in {\mathcal W}).
\end{align*}
Hence, combining the above estimates with (\ref{L^2 norms of nabla_xv in D}) and 
(\ref{emission condition}), we obtain a similar estimate 
to (\ref{estimate of the L^2 norm of the gradient from the below No1}):
\begin{align}
\V{\nabla_xv(x; \tau)}^2 \geq C\tau^2\Big\{
&\int_{{\mathcal W}\cap(\overline{D}{\times}\overline{B}{\times}\overline{B})}
e^{-\tau(l(x, y)+l(x, \xi))}dxdyd{\xi}
\label{estimate of the L^2 norm of the gradient from the below for the transmission case}
\\&\hskip20mm
- C'\V{f}_{L^2(B)}^2e^{-2{\tau}l(D, B)-c_1\tau}
\Big\}. 
\nonumber
\end{align}

\par

Since $(x^{(1)}, y^{(1)})$ attains $l(x, y) = l(D, B)$, it follows that
$x^{(1)} \in \partial{D}$ and $y^{(1)} \in \partial{B}$. 
From condition (C), there exist cones $C_{x^{(1)}}$ with vertex $x^{(1)}$ and 
$C_{y^{(1)}}$ with vertex $y^{(1)}$ such that
$C_{x^{(1)}}\setminus\{x^{(1)}\} \subset D\cap{B_\delta(x^{(1)})}$ and 
$C_{y^{(1)}}\setminus\{y^{(1)}\} \subset B\cap{B_\delta(y^{(1)})}$.

\par

We take $\tilde{z}'_1 = \tilde{z}'(x^{(1)}, y^{(1)})$. It follows that
\begin{align*}
l(x, y) &\leq 
\frac{\v{x - \tilde{z}_1'}}{\sqrt{\gamma_-}}+\frac{\v{y - \tilde{z}_1'}}{\sqrt{\gamma_+}}
\leq \frac{\v{x - x^{(1)}}+\v{x^{(1)} - \tilde{z}_1'}}{\sqrt{\gamma_-}}
+ \frac{\v{y - y^{(1)}}+\v{y^{(1)} - \tilde{z}_1'}}{\sqrt{\gamma_+}}
\\&
= \frac{\v{x^{(1)}-\tilde{z}_1'}}{\sqrt{\gamma_-}} 
+ \frac{\v{y^{(1)} - \tilde{z}_1'}}{\sqrt{\gamma_+}}
+ \frac{\v{x - x^{(1)}}}{\sqrt{\gamma_-}}
+ \frac{\v{y - y^{(1)}}}{\sqrt{\gamma_+}}
\\&
= l(x^{(1)}, y^{(1)}) + \frac{\v{x - x^{(1)}}}{\sqrt{\gamma_-}}
+ \frac{\v{y - y^{(1)}}}{\sqrt{\gamma_+}},
\end{align*}
which yields
$$
l(x, y) \leq l(D, B) + \frac{\v{x - x^{(1)}}}{\sqrt{\gamma_-}}
+ \frac{\v{y - y^{(1)}}}{\sqrt{\gamma_+}}.
$$
Noting the above estimate and 
${\mathcal W} \supset B_{\delta}(x^{(1)}){\times}B_{\delta}(y^{(1)}){\times}B_{\delta}(y^{(1)})$, 
from Lemma \ref{estimate of integrals on cones}, we obtain
\begin{align*}
\int_{{\mathcal W}\cap(\overline{D}{\times}\overline{B}{\times}\overline{B})}&
e^{-\tau(l(x, y)+l(x, \xi))}dxdyd{\xi}
\geq \int_{C_{x^{(1)}}}dx\left(\int_{C_{y^{(1)}}}e^{-{\tau}l(x, y)}dy\right)^2
\\&
\geq \int_{C_{x^{(1)}}}dx\left(\int_{C_{y^{(1)}}}
e^{-\tau\v{x - x^{(1)}}/\sqrt{\gamma_-}}e^{-{\tau}l(D, B)}
e^{-\tau\v{y - y^{(1)}}/\sqrt{\gamma_+}}dy\right)^2
\\&
= \int_{C_{x^{(1)}}}e^{-2\tau\v{x - x^{(1)}}/\sqrt{\gamma_-}}dx
\left(\int_{C_{y^{(1)}}}e^{-\tau\v{y - y^{(1)}}/\sqrt{\gamma_+}}dy\right)^2
e^{-2{\tau}l(D, B)}
\\&
\geq C\tau^{-9}e^{-2{\tau}l(D, B)}.
\end{align*}
This estimate and 
(\ref{estimate of the L^2 norm of the gradient from the below for the transmission case}) 
imply the lower estimate of Proposition \ref{estimate of nabla_xv for transmission case}. 
\hfill $\square$

\setcounter{equation}{0}
\section{Another choice of the indicator functions}
\label{Another choice of the indicator functions}

In sections \ref{no layered case} and \ref{Inverse problems in two-layered medium}, we introduce 
the indicator function $I_\tau$ defined by (\ref{Indicator function for CP}). In 
(\ref{Indicator function for CP}), $w$ is given by (\ref{def of w}). 
Note that $v$ is the solution of the reduced problem (\ref{Equation of v for CP}) or 
(\ref{Equation of v for transmission case}), which is given by the Laplace transform of 
the solution $u_0$ of the time dependent problem for the background media 
(\ref{equation with inclusion}) with $\gamma = I_3 $ 
and (\ref{equation of free case}), respectively.

\par

For another indicator function, we consider another choice of $v$. 
According to the definition of $w$, we put
$$
\tilde{v}(x; \tau) = \int_0^Te^{-{\tau}t}u_0(t, x)dt, 
$$
and introduce another indicator function $\tilde{I}_\tau$ by
\begin{equation}
\tilde{I}_\tau = \int_{0}^{T}\int_{\R^3}e^{-{\tau}t}f(x)(u(t, x)-u_0(t, x))dxdt
= \int_{\R^3}f(x)(w(x; \tau)-\tilde{v}(x; \tau))dx.
\label{another indicator function}
\end{equation}
It seems that this choice is more consistent with the definition of $w$ by
using the time integral in $[0, T]$. Moreover, $\tilde{I}_\tau$ is equivalent to 
$I_\tau$ in some sense.
\begin{prp}\label{equivalence of the indicator functions}
There exists a constant $C > 0$ such that
$$
\v{\tilde{I}_\tau - I_\tau} \leq C\tau^{-1}e^{-{\tau}T}\V{f}^2_{L^2(\R^3)}
\qquad(\tau \geq 1). 
$$
\end{prp}
\noindent
{\quad\it Proof.} 
We show Proposition \ref{equivalence of the indicator functions} for the case in 
section \ref{Inverse problems in two-layered medium}. The other case is shown similarly. 
\par
Note that $\tilde{v}$ is the weak $L^2$-solution of 
$$
(L_{\gamma_0}-\tau^2)\tilde{v}(x, \tau) + f(x) 
= e^{-{\tau}T}(\partial_tu_0(T, x)+{\tau}u_0(T, x)), 
$$
which yields that
\begin{align}
(L_{\gamma_0}-\tau^2)(\tilde{v}(x, \tau) -v(x, \tau)) = F(x; \tau)
\qquad \text{in  } \R^3  
\label{equation of the difference}
\end{align}
 in the weak sense, where $F(x; \tau) = e^{-{\tau}T}(\partial_tu_0(T, x)+{\tau}u_0(T, x))$.     
Thus (\ref{equation of the difference}) means  
$$
\int_{\R^3}\{\gamma_0\nabla_x(\tilde{v}-v)\cdot\nabla_x\phi+\tau^2(\tilde{v}-v)\phi\}dx
= \int_{\R^3}F(x; \tau)\phi(x)dx
\qquad(\phi \in H^1(\R^3)).  
$$
From the above it follows that
\begin{align*}
\int_{\R^3}\{\gamma_0\nabla_x(\tilde{v}-v)&\cdot\nabla_x(\tilde{v}-v)+\tau^2\v{\tilde{v}-v}^2\}dx
= \int_{\R^3}F(x; \tau)(\tilde{v}-v)dx
\\&
\leq \frac{1}{2\tau^2}\V{F(\cdot; \tau)}_{L^2(\R^3)}^2
+ \frac{\tau^2}{2}\V{\tilde{v}-v}_{L^2(\R^3)}^2.
\end{align*}
Hence, we have
$$
\int_{\R^3}\gamma_0\nabla_x(\tilde{v}-v)\cdot\nabla_x(\tilde{v}-v)dx
+ \frac{\tau^2}{2}\V{\tilde{v}-v}_{L^2(\R^3)}^2
\leq \frac{1}{2\tau^2}\V{F(\cdot; \tau)}_{L^2(\R^3)}^2, 
$$
which yields
\begin{align*}
&\V{\nabla_x(v - \tilde{v})}_{L^2(\R^3)}^2+\tau^2\V{v - \tilde{v}}_{L^2(\R^3)}^2
\\&\hskip30mm
\leq C\tau^{-2}e^{-2{\tau}T}(\V{\partial_tu_0(T, \cdot)}_{L^2(\R^3)}^2
+{\tau}^2\V{u_0(T, \cdot)}_{L^2(\R^3)}^2)
\end{align*}
since $\gamma_0(x) \geq c_0$ a.e. $x \in \R^3$. 
\par

From the conservation of energy for $u_0$ it follows that 
$$
\V{\partial_tu_0(T, \cdot)}_{L^2(\R^3)} \leq C\V{f}_{L^2(\R^3)}, \quad
\V{u_0(T, \cdot)}_{L^2(\R^3)} \leq CT\V{f}_{L^2(\R^3)}, \quad
$$
which yields 
\begin{align*}
\V{\nabla_x(v - \tilde{v})}_{L^2(\R^3)} + \tau\V{v - \tilde{v}}_{L^2(\R^3)}
&\leq C\V{f}_{L^2(\R^3)}e^{-{\tau}T}
\quad(\tau > 0).
\end{align*}
From this estimate it follows that
\begin{align*}
\vert \tilde{I}_\tau - I_\tau \vert \leq \Big\vert\int_{\R^3}f(x)(v - \tilde{v})dx\Big\vert
\leq \V{f}_{L^2(\R^3)}\V{v - \tilde{v}}_{L^2(\R^3)}
\leq C\V{f}_{L^2(\R^3)}^2e^{-{\tau}T}\tau^{-1},
\end{align*}
which completes the proof of Proposition \ref{equivalence of the indicator functions}.
\hfill $\square$

From Proposition \ref{equivalence of the indicator functions}, we can replace $I_\tau$ in 
(\ref{estimate of I_tau for CP -case}) and (\ref{estimate of I_tau for CP +case}) by 
$\tilde{I}_\tau$. Hence, mathematically, both indicator functions are equivalent to each other. 
It seems that adoption of $\tilde{I}_\tau$ is a consistent approach to the problem.  
However, in practice, $I_\tau$ is 
better since we can use the fundamental solution of the reduced problems directly.



%
%

\footnotesize
\par\noindent
{\sc Kawashita, Mishio}
\par\noindent
{\sc Mathematics Program, \\
Graduate School of Advanced Science and Engineering,
\\
Hiroshima University, Higashihiroshima 739-8526, Japan.}
\vskip1pc
\footnotesize
\par\noindent
{\sc Kawashita, Wakako}
\par\noindent
{\sc Electrical, Systems, and Control Engineering Program, \\
Graduate School of Advanced Science and Engineering, \\
Hiroshima University, 
Higashihiroshima 739-8527, Japan.}

\end{document}